\documentclass[12pt]{amsart} 
\usepackage[margin=1.3in]{geometry}
\usepackage{graphicx}
\usepackage{multirow}
\usepackage[table,xcdraw]{xcolor}
\usepackage{graphicx,natbib}
\usepackage{amsmath}
\usepackage{placeins}
\usepackage{amssymb}
\usepackage{epstopdf}
\usepackage{pdfpages}
\usepackage{amsthm}
\usepackage{xcolor}
\usepackage{setspace}
\usepackage[colorlinks, linkcolor=violet,citecolor=violet]{hyperref}

\newtheorem{theorem}{Theorem}[section]

\newtheorem{remark}[theorem]{Remark}

\newcommand*\Prob{\mathop{}\!\mathbb{P}}
\newcommand*\E{\mathop{}\!\mathbb{E}}

\newcommand*\dd{\mathop{}\!\mathrm{d}}

\usepackage{multicol}

\title{Trimmed extreme value estimators for censored heavy-tailed data}

\author[M. \smash{Bladt}]{Martin Bladt}
\address[M. Bladt]{Department of Actuarial Science, Faculty of Business and Economics, University of Lausanne, CH-1015 Lausanne, Switzerland}

\email{{martin.bladt@unil.ch}}

\author[H. \smash{Albrecher}]{Hansj\"org Albrecher}
\address[H. Albrecher]{ Department of Actuarial Science, Faculty of Business and Economics and Swiss Finance Institute, University of Lausanne, CH-1015 Lausanne, Switzerland}

\email{{hansjoerg.albrecher@unil.ch}}

\author[J. \smash{Beirlant}]{Jan Beirlant}
\address[J. Beirlant]{Department of Mathematics,
	KU Leuven, Celestijnenlaan 200B, B-3001 Leuven, Belgium and  Department of Mathematical Statistics and Actuarial Science, University of the Free State, South Africa}

\email{{jan.beirlant@kuleuven.be}}

\begin{document}
\maketitle
\begin{abstract}
We consider estimation of the extreme value index and extreme quantiles for heavy--tailed data that are right-censored. We study a general procedure of removing low importance observations in tail estimators.
This trimming procedure is applied to the state-of-the-art estimators for randomly right-censored tail estimators. Through an averaging procedure over the amount of trimming we derive new kernel type estimators.  Extensive simulation suggests that one of the new considered kernels leads to a highly competitive estimator against virtually any other available alternative in this framework. Moreover we propose an adaptive selection method for the amount of top data used in estimation based on the trimming procedure minimizing the asymptotic mean squared error. We also provide an illustration of this approach to simulated as well as to real-world MTPL insurance data.
\end{abstract} 
\keywords{Keywords: heavy-tailed estimation; right-censored data; asymptotic distributions}

\section{Introduction}

In recent years the problem of tail estimation with right-censored data  has received considerable attention starting with \cite{beirlcens2007} and \cite{einmahl2008statistics}, who considered the problem for different domains of attraction. Most efforts however have been dedicated to heavy-tailed distributions, with several papers being motivated by heavy-tailed insurance claim data with long development times of the claims, see e.g. \cite{beirlant10,beirlant2016bias,beirlant2018penalized,beirlant19,worms2014new,worms2016lynden,worms2018extreme,ndao2014nonparametric}.
The underlying model assumption here is that the random variable of interest $X$ has a Pareto-type distribution function
\begin{align}\label{purepareto}
F(x)=\Prob(X\le x)=1-x^{-1/\xi}\ell(x), \quad \xi>0, x>1,
\end{align}
where $\ell$ is slowly varying at infinity:
\begin{align*}
\lim_{x\to\infty}\frac{\ell(tx)}{\ell(x)}=1, \mbox{ for every } t>1.
\end{align*}
\cite{worms2019weibull} discuss the analogous problem for Weibull-type distributions. Right-censoring for heavy-tailed distributions in a regression setting was discussed in \cite{ndao2016nonparametric,dierckx19,goegebeur19a,stupfler16}, while \cite{goegebeur19b} is on a bivariate extension. \cite{stupfler19} discusses the case of dependent censoring. \\

Here we revisit the case of heavy-tailed data.   More specifically, we consider the random right-censoring model,  where the independent and identically distributed (i.i.d.)\ observations  $X_1, \dots, X_n$ of $X$ may be preceded by censoring variables $C_1, \dots, C_n $, and it is known if that happens. One then observes
\begin{align*}
Z_i=\min\{X_i,C_i\},\quad e_i=1\{X_i\le C_i\},\quad i=1,\dots,n,
\end{align*}
where $C_1,\dots, C_n$ is an i.i.d.\ sequence of censoring random variables, independent of the observations $X_i$. Here $1\{X_i\le C_i\}$ denotes the indicator of the event $\{X_i\le C_i\}$. {\color{black}In order to obtain non-degenerate and tractable identities, one assumes that also the censoring variables are Pareto-type distributed with distribution function
\begin{align*}
G(x)=\Prob(C_1\le x)=1-x^{-1/\xi_c}\ell_c(x),\quad \quad \xi_c>0, x> 1,
\end{align*}
with $\ell_c$ another slowly varying function at infinity. This choice also motivates the asymptotic Pareto behaviour to be introduced in the sequel.} Then we have that for $x>1$,
\begin{align*}
H(x)=\Prob(Z_1\le x)=1-x^{-1/\xi_z}\ell_z(x),\quad \xi_z=\frac{\xi\xi_c}{\xi+\xi_c},
\end{align*}
where $\ell_z(x)=\ell(x)\ell_c(x)$. 
As explained in Einmahl et al. (2008), the parameter $p= \xi_z/\xi= \frac{\xi_c}{\xi +\xi_c}$ is the limit of $p(z):=\Prob (e_1 =1 | Z_1=z)$ as $z \to \infty$, and can be interpreted as the non-censoring probability in the limit, or the tail limiting proportion of non-censored data. In the exact Pareto setting (i.e. $\ell$ and $\ell_c$ being constant) the censoring indicators  $e_1,\dots, e_n$ turn out to be i.i.d.\ Bernoulli($p$) random variables, independent of $Z_1,\dots, Z_n$. \\

Within this censoring and regularly varying context, \cite{beirlcens2007}  proposed a first estimator of $\xi$ in the spirit of the classical Hill estimator . Concretely, define the order statistics of the observed sample as
\begin{align*}
Z_{1,n}\le \dots\le Z_{n,n},
\end{align*}
and $e_{i,n}$ the corresponding censoring indicators, $i=1,\dots,n$. Then the Hill estimator adapted for censoring is given by
\begin{align}\label{hill}
H_k=\frac{\sum_{i=1}^k \log(Z_{n-i+1,n}/Z_{n-k,n})}{\sum_{i=1}^{k}e_{n-i+1,n}}=\frac{H_k^Z}{p_k},  \quad 1<k<n,
\end{align}
with 
\[
p_k=\frac 1k \sum_{i=1}^k e_{n-i+1,n}
\]
the fraction of non-censored data in the top $k$ observations, and 
\[
H_k^Z=\frac{1}{k}\sum_{i=1}^k \log(Z_{n-i+1,n}/Z_{n-k,n})
\]
the classical Hill estimator (Hill, 1975) based on the top $k$ observations.
\cite{einmahl2008statistics} showed that, under some regularity assumptions, $H_k$ is consistent  whatever the value of $p \in (0,1)$: 
\[
H_k=\frac{H_k^Z}{p_k}\to_p \frac{\xi_z}{p}=\xi
\]
as $k,n \to \infty$ and $k/n \to 0$. Moreover \cite{einmahl2008statistics} derived the asymptotic normality of  $H_k$ under general conditions. \\

\cite{worms2014new} proposed an alternative generalization of the Hill estimator based on the fact that 
\begin{align}\label{lmef}
\E(\log(Z/t)|Z>t)=\int_1^\infty\frac{\overline F(ut)}{\overline F(t)} \frac{1}{u}\dd u\to\xi \mbox{ as } t\to\infty,
\end{align}
where $\overline F(x)=1-F(x)$. In the exact Pareto case, the above limit is an equality. Replacing $\overline F$ with the Kaplan-Meier estimator
\begin{align}\label{KM}
\widehat{\overline F}(x)=\prod_{Z_{j,n}\le x}\left(\frac{n-j}{n-j+1}\right)^{e_{j,n}}
\end{align}
for $\overline F(x)$  yields the estimator
\begin{align}\label{worms}
H^W_k=\sum_{i=1}^k\frac{\widehat{\overline F}(Z_{n-i,n})}{\widehat{\overline F}(Z_{n-k,n})}\log(Z_{n-i+1,n}/Z_{n-i,n}),
\end{align}
which was shown to be consistent in \cite{worms2014new}.  Observe that both \eqref{hill} and \eqref{worms} reduce to $H_k^Z$ when there is no censoring. Based on simulation studies, see e.g. \cite{beirlant2018penalized}, the estimator $H_k^W$ is known to exhibit superior behaviour in comparison with $H_k$, especially with respect to bias. {\color{black} The mathematical treatment of this estimator turns out to be difficult and the asymptotic normality of $H_k^W$ has only be derived in \cite{beirlant19}  under  light censoring, i.e.  $\xi< \xi_c$ or $p>1/2$, and some regularity conditions. Here we will propose a novel family of estimators on the basis of a trimming procedure that will exhibit competitive behaviour and for which the mathematical treatment is simpler. This then also yields a generalization to the censoring case of the classical class of kernel estimators introduced in \cite{csorgo85}.}\\

Before introducing the trimming
procedure, we first propose  simplified versions of $H_k^W$, that will be more amenable for the approach in the sequel.  To this end,  note that one can write
\begin{align}\label{hkw}
H^W_k=\sum_{i=1}^k\left[\prod_{j=i+1}^k(1-1/j)^{e_{n-j+1,n}}\right]\log(Z_{n-i+1,n}/Z_{n-i,n}).
\end{align}
In the exact Pareto case one has $p(z)=p$ and 
 the term in the square bracket  has expectation
\begin{align}
\E\left[\prod_{j=i+1}^k(1-1/j)^{e_{n-j+1,n}}\right]&{\color{black}=\E\left[\prod_{j=i+1}^k(1-1/j)^{e_{j}}\right]}\nonumber\\
&=\prod_{j=i+1}^k \E\left[(1-1/j)^{e_{j}}\right]
=
\prod_{j=i+1}^k(1-p/j),\label{meankm}
\end{align}
{\color{black}where we have used the fact that the $e_{n-j+1,n}$ are i.i.d. Bernoulli random variables with the same law as the $e_j$. Notice that for Pareto-type variables, this only holds asymptotically, as $k,\,n/k\to\infty$.}
The R\'enyi representation also entails  for the exact Pareto case that
\begin{align} \label{renyi}
\E(V_j)=\frac{1}{\sum_{m=i}^k m^{-1}}\E(\log(Z_{n-j+1,n}/Z_{n-k,n}))=\xi_z
\end{align}
with $V_j = j\log {Z_{n-j+1,n}\over Z_{n-j,n}}$.
Based on \eqref{meankm} and \eqref{renyi} and using the approximations $1-p_k/j\approx \exp(-p_k/j)$ and $\sum_{j=i+1}^k j^{-1}\approx \log((k+1)/i)$, we propose the following estimator which is closely related to $H_k^W$: 
\begin{align}\label{aux2}
 H^A_k=\frac{1}{k+1} \sum_{i=1}^k\left(\frac{i}{k+1}\right)^{p_k-1}\frac{1}{\log((k+1)/i)}\log(Z_{n-i+1,n}/Z_{n-k,n}), \quad k<n,
\end{align}
where the log-spacings in the sum are all taken with respect to the same baseline order statistic $Z_{n-k,n}$. The latter will allow to apply the trimming operation of removing low importance observations in the tail estimation,  developed in \cite{abbtrim} for the non-censoring case, to the present situation with censoring.\\


In this paper, we extend the trimming method proposed in \cite{abbtrim} to the case of random right censoring, both for $H_k$ and $H_k^A$. Averaging the trimmed statistics over the amount of trimming then leads to new estimators which belong to 	a  general  family of kernel estimators comprising $H_k$ and $H_k^A$. This family turns out to be closed under the proposed averaging operation after trimming. Then we study the basic asymptotic properties of the kernel estimators in Section \ref{sec3}, even in case of heavy censoring, i.e. $p \leq 1/2$. 
Trajectories of the trimmed statistics as a function of the amount of trimming turn out to be quite flat near the optimal threshold value minimizing the mean squared error (MSE). 
Based on this, in Section \ref{sec4} -- for the first time in this setting -- an adaptive selection method for the amount of top data used in tail estimation is proposed.
 In Sections \ref{secsim} and \ref{secins}, we show through simulations and a case study from insurance that the new kernel estimators  and the threshold selection method exhibit  promising properties.   

\section{Trimmed estimators for $\xi$}

\subsection{Trimming tail estimators}

In \cite{abbtrim}, lower trimming of the classical Hill estimator was shown to be an effective strategy to obtain Hill-type plots with lower variance arising from the changes of the baseline order statistic, which aids in the visual selection of a horizontal part of the trajectory. Here, we extend this approach to the censored case, and consider lower trimming of the estimators $H_k$ and $H^A_k$, deleting the smallest $k-b$ ($b \leq k$) peaks over thresholds $Z_{n-i+1,n}/Z_{n-k,n}$, $i=b+1, \ldots,k$: trimming $H_k^Z$ as in \cite{abbtrim} one obtains a trimmed version of $H_k$
\begin{align}
H_{b,k}=\frac{1}{1+\sum_{j=b+1}^k j^{-1}}\cdot \frac{{1 \over b}\sum_{i=1}^b \log(Z_{n-i+1,n}/Z_{n-k,n})}{p_k},\quad b \le k\le n-1,
\end{align}
while for $H^A_k$ we propose 
\begin{align}\label{wormtrim2}
 H^A_{b,k}=\frac{1}{b+1} \sum_{i=1}^b\left(\frac{i}{b+1}\right)^{p_k-1}\frac{1}{\log((k+1)/i)}
\log(Z_{n-i+1,n}/Z_{n-k,n}),\quad b\le k\le n-1,
\end{align}
since, when $p_k$ is replaced by the exact value $p$, using \eqref{renyi} the expected value equals $ \frac{\xi_z}{b+1} \sum_{i=1}^b\left(\frac{i}{b+1}\right)^{p-1} \frac{\sum_{m=i}^{k}m^{-1}}{\log((k+1)/i)} \approx \xi$ for $b$ sufficiently large.  Note also that
$H_{k,k}=H_k$ and $H_{k,k}^A=H_k^A$.

\subsection{Averaging and kernels}
The above trimming procedure naturally leads  to new estimators when considering  the empirical mean of the trimmed estimators across $b=1,\ldots,k$:
\begin{align*}
\frac 1k\sum_{b=1}^k H_{b,k}, \quad
\frac 1k\sum_{b=1}^k H_{b,k}^A.
\end{align*}
For instance, in case of  $H_{b,k}^A$
this is asymptotically equivalent  to 
\begin{align}\label{aux_bar}
\overline H^A_k=\frac1k \sum_{i=1}^k\frac{1}{1-p_k}\left(\left(\frac{i}{k+1}\right)^{p_k-1}-1\right)\frac{1}{\log((k+1)/i)}\log(Z_{n-i+1,n}/Z_{n-k,n}), \; k < n,
\end{align}
as can be seen using partial summation and a simple Riemann sum approximation
${1 \over k}\sum_{b=i}^k ({b \over k+1})^{-p_k}\approx \frac{1}{1-p_k}\left(1-\left(\frac{i}{k+1}\right)^{1-p_k}\right)$ for $k \to \infty$.\\

In fact $H_k$, $H_k^A$ and $\overline H_k^A$ can all be put into a kernel framework, by defining 
\begin{align}\label{kernel}
H^{\mathcal{K}}_k=\frac1k \sum_{i=1}^k\mathcal{K}\left(\frac{i}{k+1},p_k\right)\frac{1}{\log((k+1)/i)}\log(Z_{n-i+1,n}/Z_{n-k,n}),
\end{align}
where $\mathcal{K}$ is a positive kernel function satisfying
\begin{align*}
\int_0^1\mathcal{K}(u;p)\dd u=\frac 1p, \mbox{ for all } p\in (0,1].
\end{align*}
In particular, we get
\begin{align*}
H_k=H_k^{\mathcal{K}_0}, \quad &\text{with}\quad  \mathcal{K}_0(u,p)=\frac 1p \log\left(\frac 1u\right),\\
H_k^A=H_k^{\mathcal{K}_1}, \quad &\text{with}\quad  \mathcal{K}_1(u,p)=u^{p-1},\\
\overline H_k^A=H_k^{\mathcal{K}_2}, \quad &\text{with}\quad  \mathcal{K}_2(u,p)=\frac{u^{p-1}-1}{1-p}.
\end{align*}
Note that $H_k^W$  does not fall into this framework, but its simplified version $H_k^A$ does.\\

Also notice that, when trimming any kernel estimator $H_k^{\mathcal{K}}$ to obtain
\begin{align}\label{trimmedkernel}
H^{\mathcal{K}}_{b,k}=\frac{1}{b+1} \sum_{i=1}^b\mathcal{K}\left(\frac{i}{b+1},p_k\right)\frac{1}{\log((k+1)/i)}\log(Z_{n-i+1,n}/Z_{n-k,n}),
\end{align}
the averaging operation $\frac 1k \sum_{b=1}^k H^{\mathcal{K}}_{b,k}$ leads to an associated kernel estimator 
 \begin{align}\label{averagekernel}
H^{\bar{\mathcal{K}}}_k=\frac1k \sum_{i=1}^k\bar{\mathcal{K}}\left(\frac{i}{k+1},p_k\right)\frac{1}{\log((k+1)/i)}\log(Z_{n-i+1,n}/Z_{n-k,n}),
\end{align}
with 
\begin{align*}
\bar{\mathcal{K}}(u,p)= \int_u^1 \frac{\mathcal{K}(v,p)}{v} dv, 
\end{align*}
where $\bar{\mathcal{K}}({i \over k+1},p)$ is obtained using a Riemann approximation of $\frac{1}{k}\sum_{b=i}^{k}{k \over b+1}\mathcal{K} \left({i \over k+1} {k+1 \over b+1},p\right)$ as $k \to \infty$ for fixed $i$. \\

Rewriting the kernel estimators $H_k^{\mathcal{K}}$ from \eqref{kernel} in terms of the random variables $V_j$ ($j=1,\ldots,k$) has some theoretical advantage,  since for the exact Pareto case  these  are independent and exponentially distributed with mean $\xi_z$ thanks to the R\'enyi representation:
\begin{eqnarray*}
H_k^{\mathcal{K}} &=& \frac1k \sum_{i=1}^k\mathcal{K}\left(\frac{i}{k+1},p_k\right)\frac{1}{\log((k+1)/i)}
\sum_{j=i}^k {V_j \over j} \\
&=& \frac1k \sum_{j=1}^k {V_j \over j} \sum_{i=1}^j \mathcal{K}\left(\frac{i}{k+1},p_k\right)\frac{1}{\log((k+1)/i)} \\
&=& \frac1k \sum_{j=1}^k V_j \; \widetilde{\mathcal{K}}_k \left(\frac{j}{k+1},p_k\right)
\end{eqnarray*}
with
\[
\widetilde{\mathcal{K}}_k \left(\frac{j}{k+1},p_k\right)
= {1 \over j/(k+1)} 
{1 \over k+1} \sum_{i=1}^j \mathcal{K}\left(\frac{i}{k+1},p_k\right)\frac{1}{\log((k+1)/i)}. 
\]
Using a Riemann approximation, one can also propose to use 
\begin{equation} \label{tildeK}
\widetilde{H}_k^{\mathcal{K}}=
 \frac1k \sum_{j=1}^k \widetilde{\mathcal{K}}\left(\frac{j}{k+1},p_k\right) \; V_j
\end{equation}
with associated kernel
\[
\widetilde{\mathcal{K}}\left( u,p \right)= {1 \over u}
\int_0^u {\mathcal{K}(v,p) \over \log (1/v)}dv
\]
that also satisfies the norming $\int_0^1 \widetilde{\mathcal{K}}\left( u,p \right)du =1/p$. 
The class of estimators $\widetilde{H}_k^{\mathcal{K}}$ can be considered as generalizations of the kernel estimators proposed in \cite{csorgo85} from the non-censoring to the censoring case. 


\subsection{Quantile estimation}
Following the approach from \cite{weissman1978estimation} it is possible to construct quantile estimators as a function of the sample size as follows. Let the quantile function of a regularly varying tail be $Q(p)$. The regular variation property implies that the ratio of increasingly large quantiles satisfies
\begin{align}\label{weismannasympt}
\frac{Q(1-p)}{Q(1-k/n)}\sim\left(\frac{k}{np}\right)^{\xi}, \quad \text{as}\quad p\downarrow 0, \:k/n \to 0,\: np=o(k),
\end{align}
as discussed in Section 4.3 in de Haan and Ferreira (2007). This then leads to a quantile estimator based on $k$ order statistics and the kernel $\mathcal{K}$ as \begin{align}
\hat Q_k(1-p)=\hat Q^{KM}(1-k/n)\cdot\left(\frac{k}{np}\right)^{H_k^{\mathcal{K}}},
\end{align}
where $\hat Q^{KM}$ is the quantile function derived from the Kaplan-Meier estimator $\widehat{\overline F}$ defined in \eqref{KM}.


\section{Asymptotic representations}\label{sec3}

In this section we derive the asymptotic distributions of the kernel estimators and their trimmed counterparts as introduced in the preceding section. In \cite{einmahl2008statistics} the asymptotics for $H_k = H_k^{\mathcal{K}_0}$ was discussed in detail. \cite{beirlant19} provided an asymptotic normality result for $H_k^W$ when $p>1/2$, but that estimator is not in the current kernel framework. Here we provide asymptotic representations for the class of kernel estimators in the form $\widetilde{H}_k^{\mathcal{K}}$. To this end, we make use of second-order assumptions which were first proposed in Hall and Welsh (1985) and which have widely been used in papers on the estimation of the extreme value index for Pareto-type  distributions both in the non-censoring case such as  \cite{csorgo85} and the censoring case  as in \cite{beirlant19}:
\begin{align} \label{HW}
\ell(x)=C(1+Dx^{-\beta}(1+o(1)), \;\;
\ell_c(x)&=C_c(1+D_cx^{-\beta_c}(1+o(1)),\; x \to \infty,
\end{align}
where $\beta,\beta_c,C,C_c$ are positive constants and $D,D_c$ are real constants. It now follows that
\begin{align*}
\ell_z(x)=C_z(1+D_zx^{-\beta_z}(1+o(1)),
\end{align*}
where
\begin{align*}
C_z=CC_c,\quad \beta_z=\min\{\beta,\beta_c\}, \quad D_z=D\cdot1_{\beta\leq\beta_c}+D_c\cdot1_{\beta_c\leq\beta}.
\end{align*}
Further we set $\rho_z=-\beta_z\xi_z$, $Q_{0,z}(t)=-\xi_z^2 \beta_z D_z C^{\rho_z}t^{\rho_z}$.
Concerning the scaled spacings  $V_j, \; j=1,\ldots,k$, one then has   the following expansion as $n,k \to \infty$ and $k/n \to 0$ as given in Theorem 4.1 in \cite{bgst}:
\begin{equation}\label{V_exp}
V_j = \left( \xi_z +Q_{0,z}(n/k)(\frac{j}{k+1})^{-\rho_z}\right)E_j + R_{j,n} +o_p(Q_{0,z}(n/k)),
\end{equation}
with $E_j$ standard exponential random variables, independent with each $n$, and $\left|\sum_{j=i}^k R_{j,n}/j \right|=
o_p(Q_{0,z}(n/k))\max (\log ((k+1)/i,1))$.  
Next, from \cite{einmahl2008statistics} and \cite{beirlant2016bias}  one obtains that 
\begin{align}\label{pk_expansion}
\sqrt{k} (p_k-p) =& \sqrt{p(1-p)}N(1+o_p(1))  + \sqrt{k}Q_{0,z}(n/k)\frac{\kappa_z}{1-\rho_z}(1+o_p(1)),
\end{align}
{\color{black}where $N \sim N(0,1)$ can be chosen appropriately} independent of $\{V_1,V_2,\ldots \}$ , and $\kappa_z= -\frac{(D\xi)_z}{D_z\xi \xi_c}$, with \\
$(D\xi)_z = (D\xi)1_{\beta \leq \beta_c}- (D_c\xi_c)1_{\beta_c \leq \beta}$.
Based on \eqref{V_exp} and \eqref{pk_expansion}  we now derive that
\begin{eqnarray} \label{a1}
 \widetilde{H}_{k}^{\mathcal{K}}-\xi  &=&
{1 \over k}\sum_{j=1}^k \left( \widetilde{\mathcal{K}}({j \over b+1},p_k)-\widetilde{\mathcal{K}}({j \over k+1},p) \right) \; V_j
\nonumber \\
&& + \left( {1 \over k}\sum_{j=1}^k  \widetilde{\mathcal{K}}({j \over k+1},p ) \;  V_j
-\xi \right) \nonumber \\
&=:& T_{1,k}+ T_{2,k}. 
\end{eqnarray}
Using the mean value theorem, we have from  \eqref{pk_expansion} that
\begin{eqnarray}
T_{1,k} &\sim_p& \xi_z  (p_k-p) \, \alpha_b^{\widetilde{\mathcal{K}}} \nonumber \\ &=& \xi_z \alpha_b^{\widetilde{\mathcal{K}}}
\left( \sqrt{p(1-p)}{N \over \sqrt{k}} +Q_{0,z}(n/k)\frac{\kappa_z}{1-\rho_z}(1+o(1))\right),\label{a2}
\end{eqnarray}
with  $\alpha_k^{\widetilde{\mathcal{K}}} = {1 \over k} \sum_{j=1}^k {\partial \widetilde{\mathcal{K}} \over \partial p}({j \over k+1},p)$.
Next, using \eqref{V_exp},
\begin{eqnarray}
T_{2,k} &=& \xi_z \; \left( {1 \over k}\sum_{j=1}^k \widetilde{\mathcal{K}}_k ({j \over k+1},p) - {1 \over p} \right) \nonumber \\
&& + \xi_z  \; {1 \over k}\sum_{j=1}^k  \widetilde{\mathcal{K}}_k({j \over k+1},p) (E_j-1)\nonumber \\
&& +  Q_{0,z}(n/k)  \; {1 \over k}\sum_{j=1}^k  \widetilde{\mathcal{K}}_k ({j \over k+1},p)\left( {j \over k+1}\right)^{-\rho_z} E_j \nonumber \\
&&+ {1 \over k}\sum_{j=1}^k R_{j,n} \widetilde{\mathcal{K}}_k ({j \over k+1},p) + o_p(Q_{0,z}(n/k)). \label{a3}
\end{eqnarray}
Concerning the second last term in \eqref{a3} we find that
\begin{eqnarray*}
{1 \over k}\sum_{j=1}^k R_{j,n} \widetilde{\mathcal{K}}_k ({j \over k+1},p) &=& {1 \over k}\sum_{j=1}^k {R_{j,n} \over j} \sum_{i=1}^j \frac{\mathcal{K}({i \over k+1},p)}{\log ((k+1)/i)} \\
&=& {1 \over k}\sum_{i=1}^k \frac{\mathcal{K}({i \over k+1},p)}{\log ((k+1)/i)} \sum_{j=i}^k {R_{j,n} \over j} =o_p(Q_{0,z}(n/k)).
\end{eqnarray*}

From \eqref{a2} and \eqref{a3} we can now state an asymptotic expansion for $H_k^{\mathcal{K}} - \xi$.
\begin{theorem} Under \eqref{HW} we have as $k,n \to \infty$ and $k/n \to 0$
\begin{align}  
H_k^{\mathcal{K}} - \xi  
=& -\sqrt{p(1-p)}{\xi_z \over p^2}{N \over \sqrt{k}}  
+  \xi_z  \; {1 \over k}\sum_{j=1}^k  \widetilde{\mathcal{K}}_k({j \over k+1},p) (E_j-1)
\nonumber \\
 & + Q_{0,z}(n/k) 
\left\{ {-\kappa_z\xi_z \over p^2(1-\rho_z)}
+ \int_0^1 u^{-\rho_z}\widetilde{\mathcal{K}}(u,p)du \right\}(1+o(1)) \nonumber \\
& + o_p(Q_{0,z}(n/k)) + \xi_z \; \left( {1 \over k}\sum_{j=1}^k \widetilde{\mathcal{K}}_k ({j \over k+1},p) - {1 \over p} \right),\label{expansion}
\end{align}
{\color{black}where $N$ and $\{E_j, j \geq 1\}$ are introduced respectively in \eqref{pk_expansion} and \eqref{V_exp}.}
\end{theorem}
In order to select an optimal $k$, we minimize then the following  asymptotic mean squared error  of $H_k^{\mathcal{K}}$:
\begin{align}
\label{AMSE_K}
\text{AMSE} (H_k^{\mathcal{K}})= \xi_z^2  \;v_{k,p} + 
Q_{0,z}^2 ({n \over k})\; b_{k,p}
\end{align}
with
\begin{eqnarray*}
v_{k,p} 
&=&
{1 \over k}{1-p \over p^3}+ {1 \over k^2}
\sum_{j=1}^k  \widetilde{\mathcal{K}}^2_k({j \over k+1},p), 
\\
b_{k,p} 
&=&
\left\{ {-\kappa_z\xi_z \over p^2(1-\rho_z)}
+ \int_0^1 u^{-\rho_z}\widetilde{\mathcal{K}}(u,p)du 
\right\}^2.
\end{eqnarray*}

\begin{remark}\normalfont
{\color{black} Note that the first two terms in \eqref{expansion} concern the estimation of $p$ and $\xi$, respectively. The third term can be regarded as a bias term arising from the second order assumption, which commonly appears in classical extreme value theory. The fourth and final term is proportional to $\xi_z$ and can be regarded as a discretization error term, since the sum inside the parenthesis is a Riemann approximation to $\int_0^1 \widetilde{\mathcal{K}}_k\left( u,p \right)du =1/p$. In general, this error is small but non-zero.}

\end{remark}

\noindent
{\bf Examples.}
Let us consider the expressions of $v_{k,p}$ and $b_{k,p}$ for some of the kernels $\mathcal{K}$ considered above. 
\\
1. For the adapted Hill estimator $H_k$ with $\mathcal{K}_0 (u,p)= \log(1/u)/p$ we have that $\widetilde{\mathcal{K}}_0(u,p) = 1/p$ and {\color{black}${1 \over k^2}
\sum_{j=1}^k  \widetilde{\mathcal{K}}^2_{0,k}({j \over k+1},p) \sim k^{-1}\int_0^1  \widetilde{\mathcal{K}}^2_{0} (v,p)dv=(kp^2)^{-1}$} as $k \to \infty$ for all $p \in (0,1)$ so that $v_{k,p}={1 \over k}p^{-3}$. Moreover $ \int_0^1 u^{-\rho_z}\widetilde{\mathcal{K}}_0(u,p)du = {1 \over p}(1-\rho_z)^{-1}$.
In case $\beta \leq \beta_c$, i.e.\ when the bias is largest compared with the classical Hill estimator in case of no censoring, we have that $ -\kappa_z= \xi_c^{-1}$ and then $b_{k,p}= p^{-4}(1-\rho_z)^{-2}$. Hence, under $\beta \leq \beta_c$,
\begin{align}\label{eq:amsecensored}
\text{AMSE}(H_k)=\frac{1}{p^4}\left(p \frac{\xi_z^2}{k}+\frac{Q_{0,z}^2(n/k)}{(1-\rho_z)^2}\right). 
\end{align}

\noindent
2. The estimator $\overline H_k^A=H_k^{\mathcal{K}_2}$ with $\mathcal{K}_2(u,p)=\frac{u^{p-1}-1}{1-p}$ exhibits excellent finite sample behaviour as shown below through the simulations. Here $\widetilde{\mathcal{K}}_2 (u,p) \sim {1 \over 1-p} \frac{p^{-1}u^{p-1}-1}{\log (1/u)}$ as $u \to 0$.
\\
 For the bias we have $\int_0^1 u^{-\rho_z}\widetilde{\mathcal{K}}^2_{2}(u,p)du = \int_0^1 k_{1-p}(1/v) k_{\rho_z}(1/v){dv \over \log (1/v)} $ with $k_a (w) = (w^a -1)/a$ for any $a>0$. As $k \to \infty$
  $$b_{k,p}\sim b_{2,p}:= \left\{  {-\kappa_z\xi_z \over p^2(1-\rho_z)}
+ \int_0^1 k_{1-p}(1/v) k_{\rho_z}(1/v){dv \over \log (1/v)} \right\}^2.$$
 
\noindent
{\color{black} Only in case of weak censoring, i.e.  $p \geq 1/2$, $\int_0^1 \widetilde{\mathcal{K}}^2_{2}(u,p)du < \infty$ as $k \to \infty$, so that then} 
${1 \over k^2}
\sum_{j=1}^k  \widetilde{\mathcal{K}}^2_{2,k}({j \over k+1},p) \sim k^{-1}\int_0^1 \widetilde{\mathcal{K}}^2_{2}(u,p)du$. Hence as $\sqrt{k}Q_{0,z}(n/k) \to \lambda \in \mathbb{R}$ 
\[
\sqrt{k}(\overline H_k^A - \xi) \to_d \mathcal{N}(\lambda b_{2,p} \, ,\sigma^2_{2,p}),
\]  
with $\sigma^2_{2,p}=(1-p)/p^3 + \int_0^1 \widetilde{\mathcal{K}}^2_{2}(u,p)du$.\\

\noindent
Under heavy censoring,  i.e.  $p < 1/2$, $$\lim_{k \to \infty}k^{2p} (\log k)^2 \left( {1 \over k^2}
\sum_{j=1}^k  \widetilde{\mathcal{K}}^2_{2,k}({j \over k+1},p) \right) = \sigma^2_{h,2,p}< \infty$$ so that, as $k^{p}(\log k) Q_{0,z}(n/k) \to \lambda \in \mathbb{R}$,
\[
k^p (\log k)(\overline H_k^A - \xi) \to_d \mathcal{N}(\lambda b_{2,p} \, , \sigma^2_{h,2,p}).
\]  

\section{Optimal choice of $k$ when estimating $\xi$}\label{sec4}

Denoting the trimmed version of the Hill estimator in the fully observed case by  $H^Z_{b,k} = H_{b,k} \; p_k$, it was shown in \cite{abbtrim} that the value $k_{\text{opt}}(H^Z_k)$ of $k$ minimizing the asymptotic MSE of $H_k^Z$ satisfies 
\begin{align*}
k_{\text{opt}}(H^Z_k)=\left(\frac{K}{(1-\rho_z)^2f(\rho_z)}\right)^{\frac{-1}{1-2\rho_z}}k_{\text{opt}}(H^Z_{b,k}),
\end{align*}
for a universal constant $K$ and a specific function $f$. Here $k_{\text{opt}}(H^Z_{b,k})$ is the optimal sample fraction  minimizing  the expectation of the  empirical variance $S_k^2 = {1 \over k}\sum_{b=1}^k \left( H_{b,k}^Z-    
{1 \over k}\sum_{b=1}^k H_{b,k}^Z
\right)^2$.
\\
On the other hand, from \eqref{eq:amsecensored} we obtain that under $\beta \leq \beta_c$ that
\begin{align} \label{koptH}
k_{\text{opt}}(H_k)=\left(\frac{p^{-1}\,K}{(1-\rho_z)^2f(\rho_z)}\right)^{-\frac{1}{1-2\rho_z}}k_{\text{opt}}(H^Z_{b,k}).
\end{align}
This means that the optimal $k$ for the estimator $H_k$ with respect to minimization of the AMSE is linked to the optimal $k$ of its trimmed versions for the minimization of the expected empirical variance in the non-censored case. A consequence of the above formula is that 
\begin{align*}
k_{\text{opt}}(H_k)=p^{\frac{1}{1-2\rho_z}}\,k_{\text{opt}}(H^Z_{k}).
\end{align*}
That is, a larger percentage of censoring leads to a higher threshold, when compared to the non-censored case. This can already be seen from the expression of the AMSE given in \eqref{eq:amsecensored}, where a smaller $p$ leads to more weight being given to the bias term. From an intuitive point of view, when dealing with censored datasets, two sources of bias have to be accounted for, and hence a smaller sample fraction $k$ is needed to control them.

In practice, for a given sample, one finds an estimate $\hat k_0=\hat{k}_{\text{opt}}(H_{b,k}^Z)$  of $k_{\text{opt}}(H^Z_{b,k})$ through minimization of $S_k^2$ over $k$, from which an adaptive choice of $k$ is found through 
\begin{align}\label{kopt_hill}
\left(\frac{p_{\hat k_0}^{-1}\,K}{(1-\hat\rho_z)^2f(\hat\rho_z)}\right)^{-\frac{1}{1-2\hat\rho_z}}\hat k_0,
\end{align}
using an estimate $\hat\rho_z$ of $\rho_z$ and replacing $p$ by $p_{\hat k_0}$. Estimators of the second-order parameter $\rho_z$ have been proposed for instance in \cite{fraga03}. Estimators exhibit a high variability and many authors consider the use of a fixed value for $\rho_z$ such as $\rho_z=-1$. In the next section we use the choices $\rho_z=-1,-1/2,-3/2$, but the results are not very sensitive to this parameter, and we here also propose to stick to the choice  $\rho_z=-1$.

\begin{remark}\normalfont
For kernels different from $\mathcal{K}_0$, if we restrict to the case $\beta \leq \beta_c$ and $p>1/2$, and consider the expressions $kv_{k,p}=: \tilde{v}_{k,p}$ and $b_{k,p}$ taken from \eqref{AMSE_K} as constant in $k$ (i.e. converging fast enough to its limit for $k \to \infty$) we find from \eqref{AMSE_K} that the optimal $k_{\text{opt}}$ for any kernel equals
$$ \left( -{1 \over 2\rho_z} \xi_z^2 /M^2 \right)^{1/(1-2\rho_z)}  n^{-2\rho_z /(1-\rho_z)}   \left( {\tilde{v}_{k,p} \over b_{k,p} }\right)^{1/(1-\rho_z)}= k_{\text{opt}}(H_k) \left( {\tilde{v}_{k,p} \over b_{k,p} (1-\rho_z)^2} \right)^{1/(1-\rho_z)},
$$
with $M = -\xi_z^2 \beta_z D_z C^{\rho_z}$, from which one can deduce that
$$k_{\text{opt}}^{\mathcal{K}} = \left( {\tilde{v}_{k,p} \over b_{k,p} }\right)^{1/(1-\rho_z)} \left( {K \over f(\rho_z)p_{\hat{k}_0} }\right) ^{-1/(1-\rho_z)}  \hat{k}_0 . $$

 \noindent This is then basically the same formula as for $H_k$, \eqref{kopt_hill}, but one has to calculate $\tilde{v}_{k,p}$ and $b_{k,p}$ at every $k$ with $p$ estimated by $p_{\hat{k}_0}$. In heavy censoring cases the rates of convergence of the estimators are different and the approach becomes  more involved. This makes the procedure significantly more computer-intensive. In practical applications, however, the difference between the optima across kernels and heavy and light censoring cases is rather small and does not justify the extra calculations.

\end{remark}

\section{Simulations}\label{secsim}

We performed simulations using the following distributions.
\begin{itemize}
\item Burr distribution with survival function
 $1-F(x) = \left( {\theta \over \theta + x^\beta }\right)^{\lambda}$ with  $(\theta, \beta,\lambda)$  taken as (10,2,1) for $X$ and (10,3,1) for $C$ so that $p<1/2$, next to (10,3,1) for $X$ and (10,2,1) for $C$ so that $p>1/2$, and (10,2,1) for both $X$ and $C$ with $p=1/2$.
 \item Fr\'echet distribution with $F(x) = \exp (- x^{-1/\xi})$ with $\xi$ taken as 1/2 for $X$ and 1/4 for $C$ and correspondingly $p<1/2$, as 1/4 for $X$ and 1/2 for $C$ and correspondingly $p>1/2$, and finally as $\xi=1/4$ for both $X$ and $C$, so that $p=1/2$.
 \item Log-gamma distribution with density $f(x) = {\lambda^\alpha \over \Gamma (\alpha)} (\log x)^{\alpha -1}x^{-\lambda -1}$ with $(\alpha,\lambda)$ taken as (3/2,2) for $X$ and (3/2,4) for $C$ so that $p<1/2$, as 
 (3/2,4) for $X$ and (3/2,2) for $C$ so that $p>1/2$, and as (3/2,4) for $X$ and $C$ so that $p=1/2$. Note that in this case the conditions \eqref{HW} are not satisfied. 
\end{itemize}
The results are based on $200$ simulations of sample size $n=200$ each. {\color{black} The first-order tail-determining parameters were chosen such that $1/\xi=3/2,2,4$, which seem to be realistic magnitudes for insurance applications, as will be illustrated in the next section. The remaining parameters were chosen in order to satisfy inequalities such as $p>1/2$ or $p\le 1/2$ but are otherwise arbitrary.}

\begin{figure}[hh]
\centering
\includegraphics[width=14cm,trim=0cm 0cm 0cm 0cm,clip]{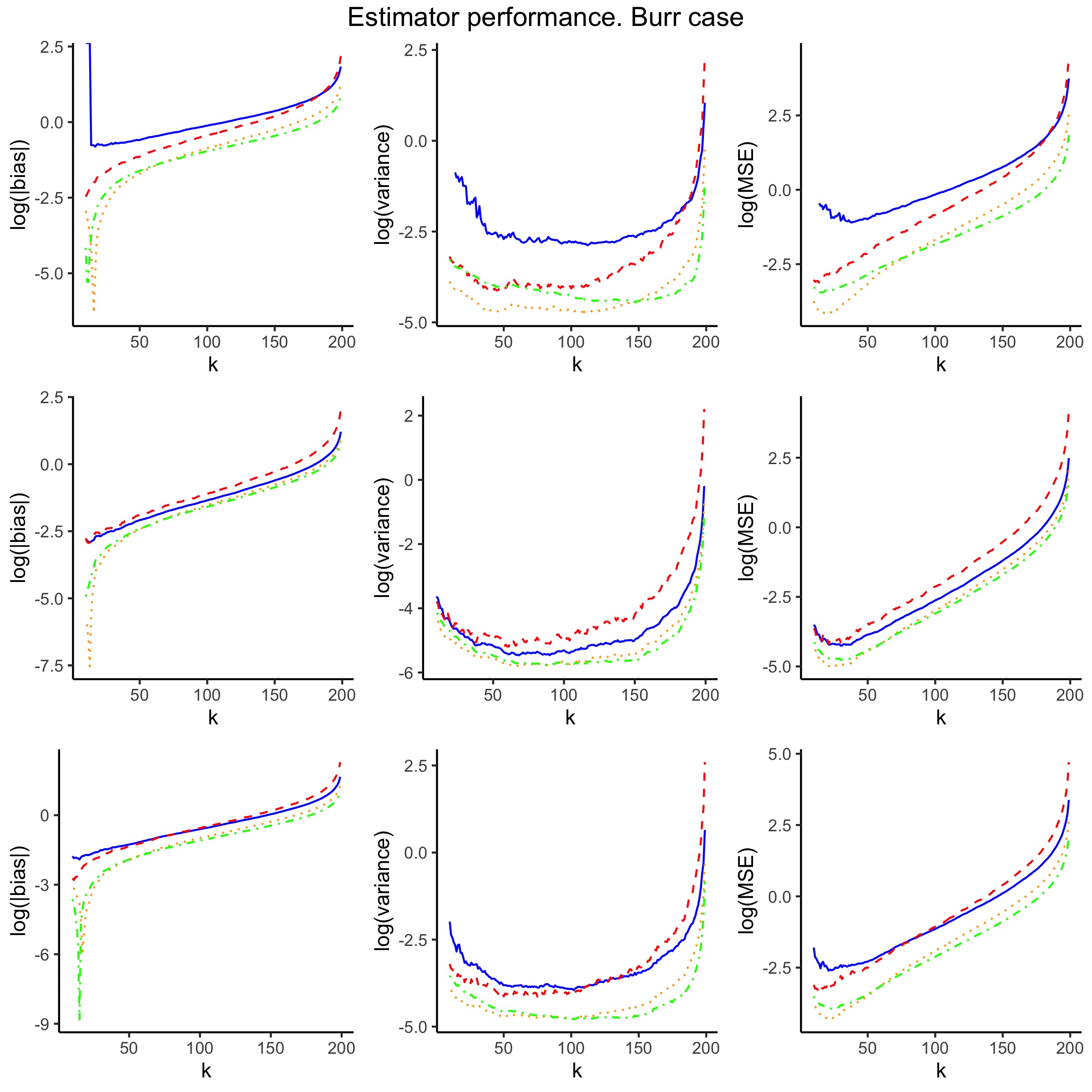}
\caption{Burr distributions: bias, variance and mean squared error of the kernel estimator ($H_k=H_k^{\mathcal{K}_0}$ in solid blue, $H^A_k=H_k^{\mathcal{K}_1}$ in dashed red, $H_k^{\mathcal{K}_2}$ in dotted orange) and the Worms estimator  $H_k^W$ (dashed and dotted green), as a function of $k$. The top, middle and bottom levels correspond to $2p<1$, $2p>1$ and $2p=1$, respectively.} 
\label{bvm_burrs}
\end{figure}

\begin{figure}[hh]
\centering
\includegraphics[width=14cm,trim=0cm 0cm 0cm 0cm,clip]{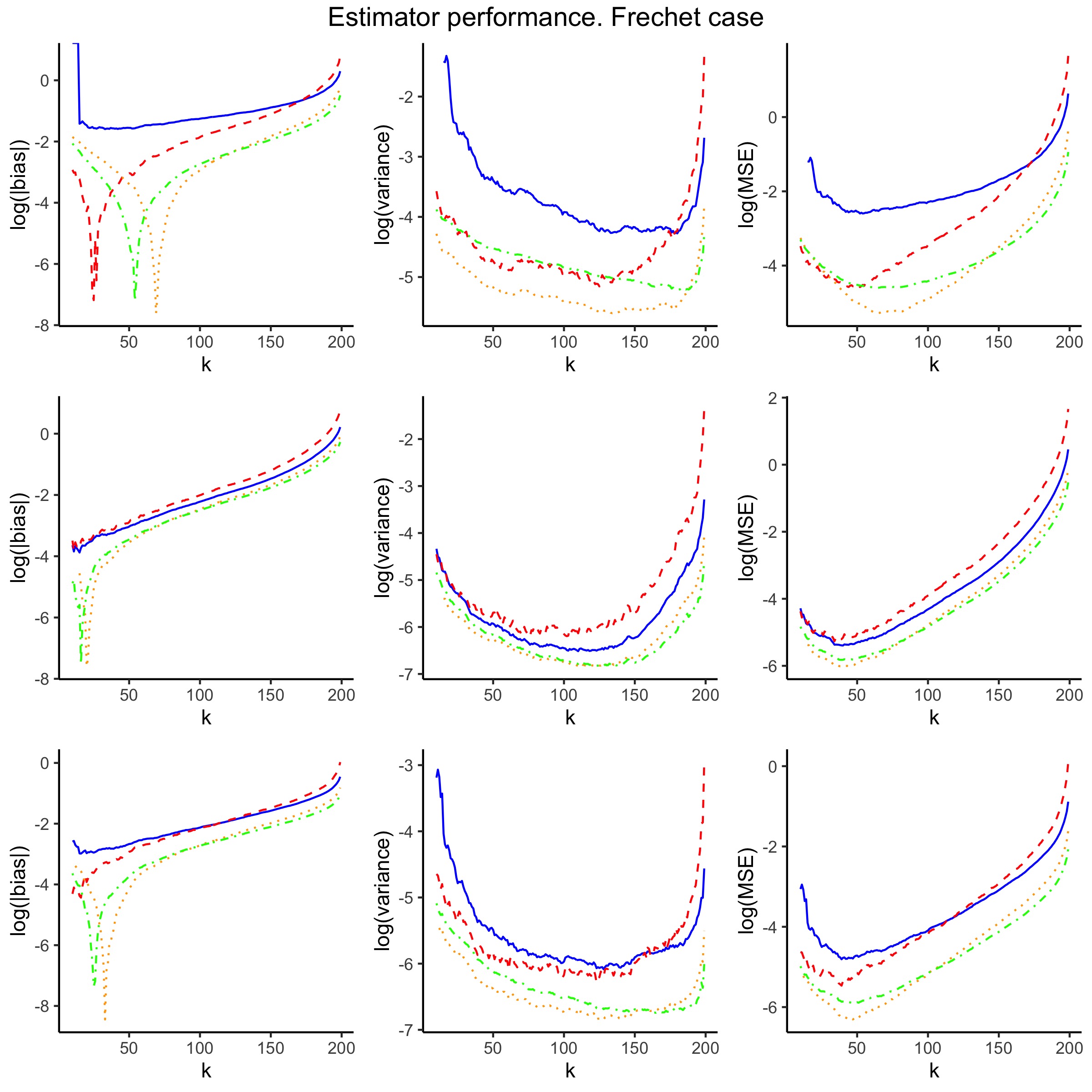}
\caption{Fr\'echet distributions: bias, variance and mean squared error of the kernel estimator ($H_k=H_k^{\mathcal{K}_0}$ in solid blue, $H^A_k=H_k^{\mathcal{K}_1}$ in dashed red, $H_k^{\mathcal{K}_2}$ in dotted orange) and the Worms estimator $H_K^W$ (dashed and dotted green), as a function of $k$. The top, middle and bottom levels correspond to $2p<1$, $2p>1$ and $2p=1$, respectively.}
\label{bvm_frechets}
\end{figure}

\begin{figure}[hh]
\centering
\includegraphics[width=14cm,trim=0cm 0cm 0cm 0cm,clip]{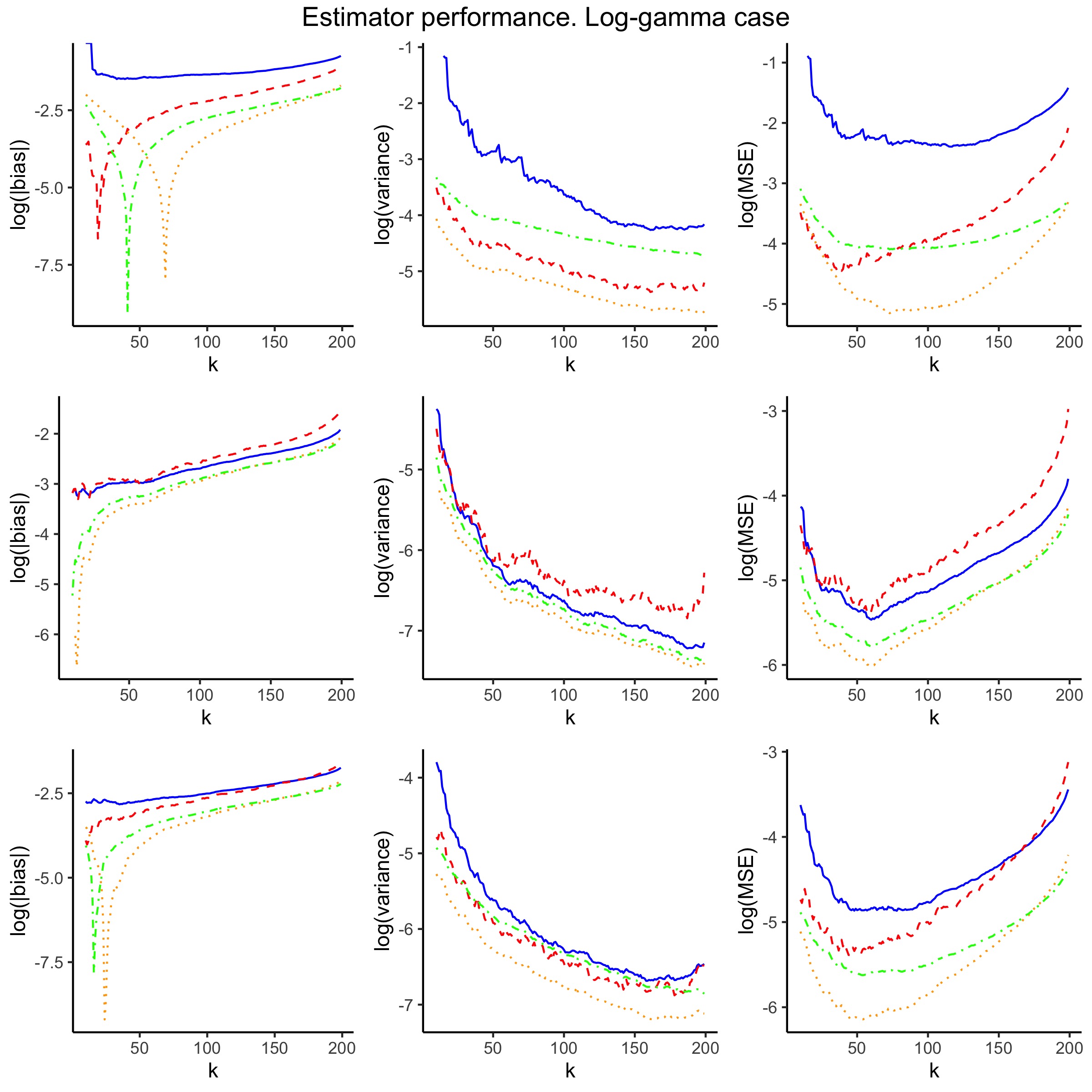}
\caption{Log-gamma distributions: bias, variance and mean squared error of the kernel estimator ($H_k=H_k^{\mathcal{K}_0}$ in solid blue, $H^A_k=H_k^{\mathcal{K}_1}$ in dashed red, $H_k^{\mathcal{K}_2}$ in dotted orange) and the Worms estimator  $K_k^W$ (dashed and dotted green), as a function of $k$. The top, middle and bottom levels correspond to $2p<1$, $2p>1$ and $2p=1$, respectively.}
\label{bvm_LGs}
\end{figure}

\begin{figure}[hh]
\centering
\includegraphics[width=14cm,trim=0cm 0cm 0cm 0cm,clip]{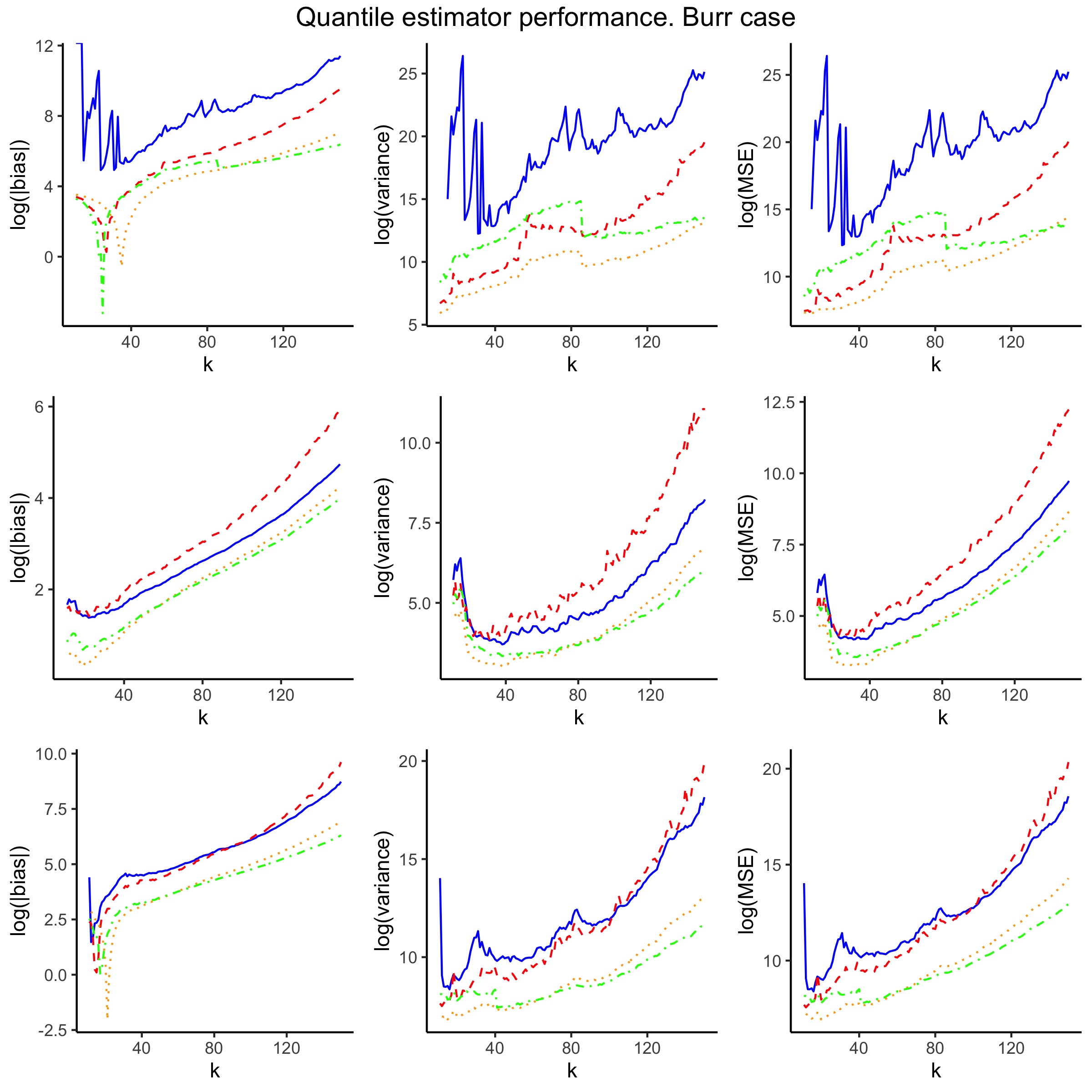}
\caption{Burr distributions. Bias, variance and mean squared error of the quantile estimator based on: the kernel estimator ($H_k=H_k^{\mathcal{K}_0}$ in solid blue, $H^A_k=H_k^{\mathcal{K}_1}$ in dashed red, $H_k^{\mathcal{K}_2}$ in dotted orange) and the Worms estimator  $H_k^W$ (dashed and dotted green), as a function of $k$. The top, middle and bottom levels correspond to $2p<1$, $2p>1$ and $2p=1$, respectively.} 
\label{VaR_bvm_burrs}
\end{figure}

\begin{figure}[hh]
\centering
\includegraphics[width=14cm,trim=0cm 0cm 0cm 0cm,clip]{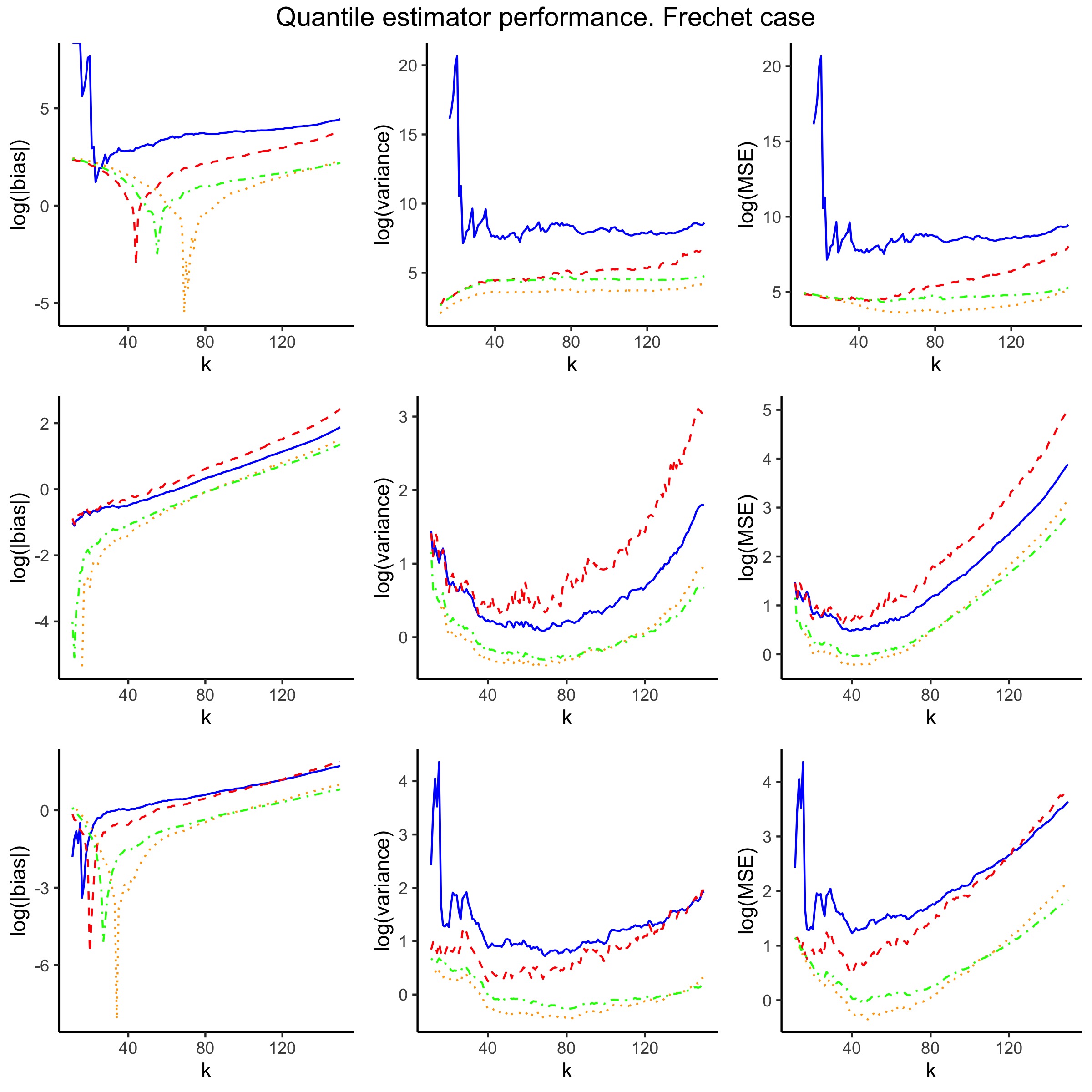}
\caption{Fr\'echet distributions. Bias, variance and mean squared error of the quantile estimator based on: the kernel estimator ($H_k=H_k^{\mathcal{K}_0}$ in solid blue, $H^A_k=H_k^{\mathcal{K}_1}$ in dashed red, $H_k^{\mathcal{K}_2}$ in dotted orange) and the Worms estimator $H_K^W$ (dashed and dotted green), as a function of $k$. The top, middle and bottom levels correspond to $2p<1$, $2p>1$ and $2p=1$, respectively.}
\label{VaR_bvm_frechets}
\end{figure}

\begin{figure}[hh]
\centering
\includegraphics[width=14cm,trim=0cm 0cm 0cm 0cm,clip]{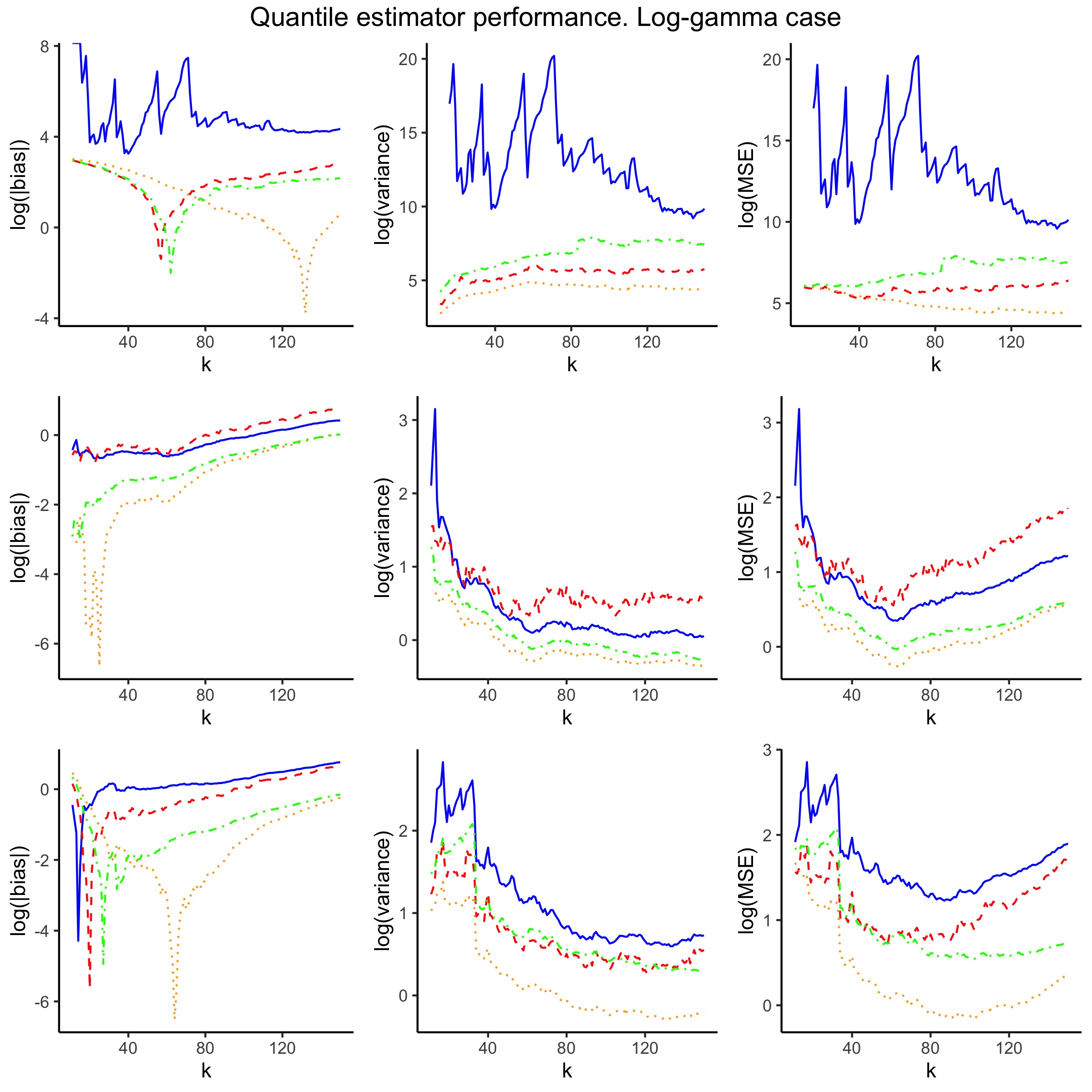}
\caption{Log-gamma distributions. Bias, variance and mean squared error of the quantile estimator based on: the kernel estimator ($H_k=H_k^{\mathcal{K}_0}$ in solid blue, $H^A_k=H_k^{\mathcal{K}_1}$ in dashed red, $H_k^{\mathcal{K}_2}$ in dotted orange) and the Worms estimator  $K_k^W$ (dashed and dotted green), as a function of $k$. The top, middle and bottom levels correspond to $2p<1$, $2p>1$ and $2p=1$, respectively.}
\label{VaR_bvm_LGs}
\end{figure}

In Figures \ref{bvm_burrs}, \ref{bvm_frechets} and \ref{bvm_LGs} we plot the bias, variance and mean squared error as a function of $k$ of the various estimators considered above. {\color{black}Observe that the plots are in logarithmic scale for display purposes. For the bias term this means that there is an asymptote corresponding to $\lim_{t\downarrow 0}\log(t)$.}
Note that the MSE characteristics of the estimator $H_k^{\mathcal{K}_2}$ are  quite comparable to those of $H_k^W$ in the Burr and Fr\'echet cases, and are even better for the log-gamma model. The corresponding analysis for the quantile estimators is given in Figures \ref{VaR_bvm_burrs}, \ref{VaR_bvm_frechets} and \ref{VaR_bvm_LGs}, which are in agreement with the former plots.

\begin{figure}[hh]
\centering
\includegraphics[width=14cm,trim=0cm 0cm 0cm 0cm,clip]{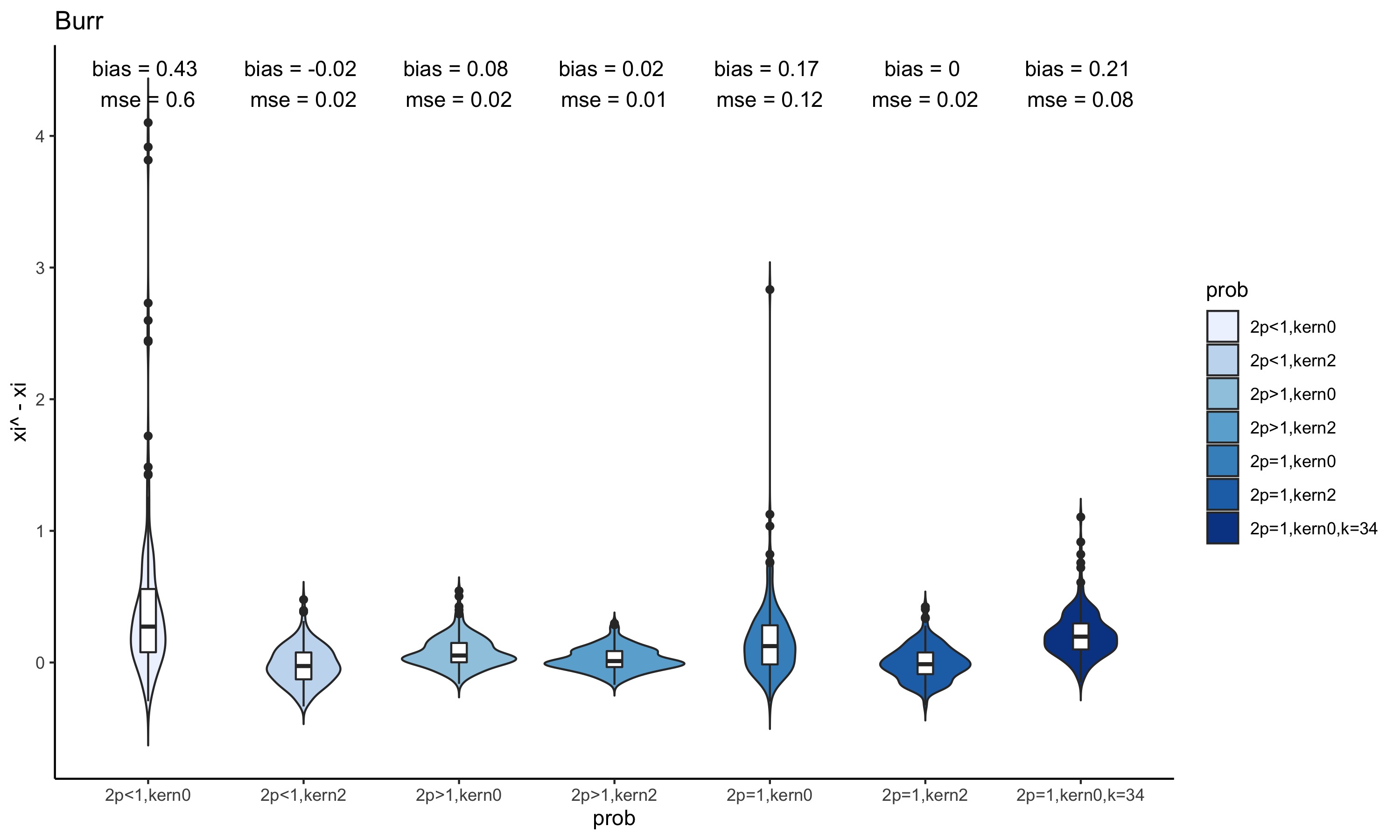}
\caption{\color{black}Violin plots for the simulation results in case of Burr distributions  under different non-censoring asymptotic probabilities. We take the difference between estimated and theoretical $\xi$. The cases $2p<1$ and $2p\ge 1$ correspond to heavy and light censoring, respectively. The labels kern0 and kern2 correspond to the use of the estimators $H_k^{\mathcal{K}_0}$ and $H_k^{\mathcal{K}_2}$. The specified $k=34$ is the theoretical optimal sample fraction.}
\label{v1}
\end{figure}

\begin{figure}[hh]
\centering
\includegraphics[width=14cm,trim=0cm 0cm 0cm 0cm,clip]{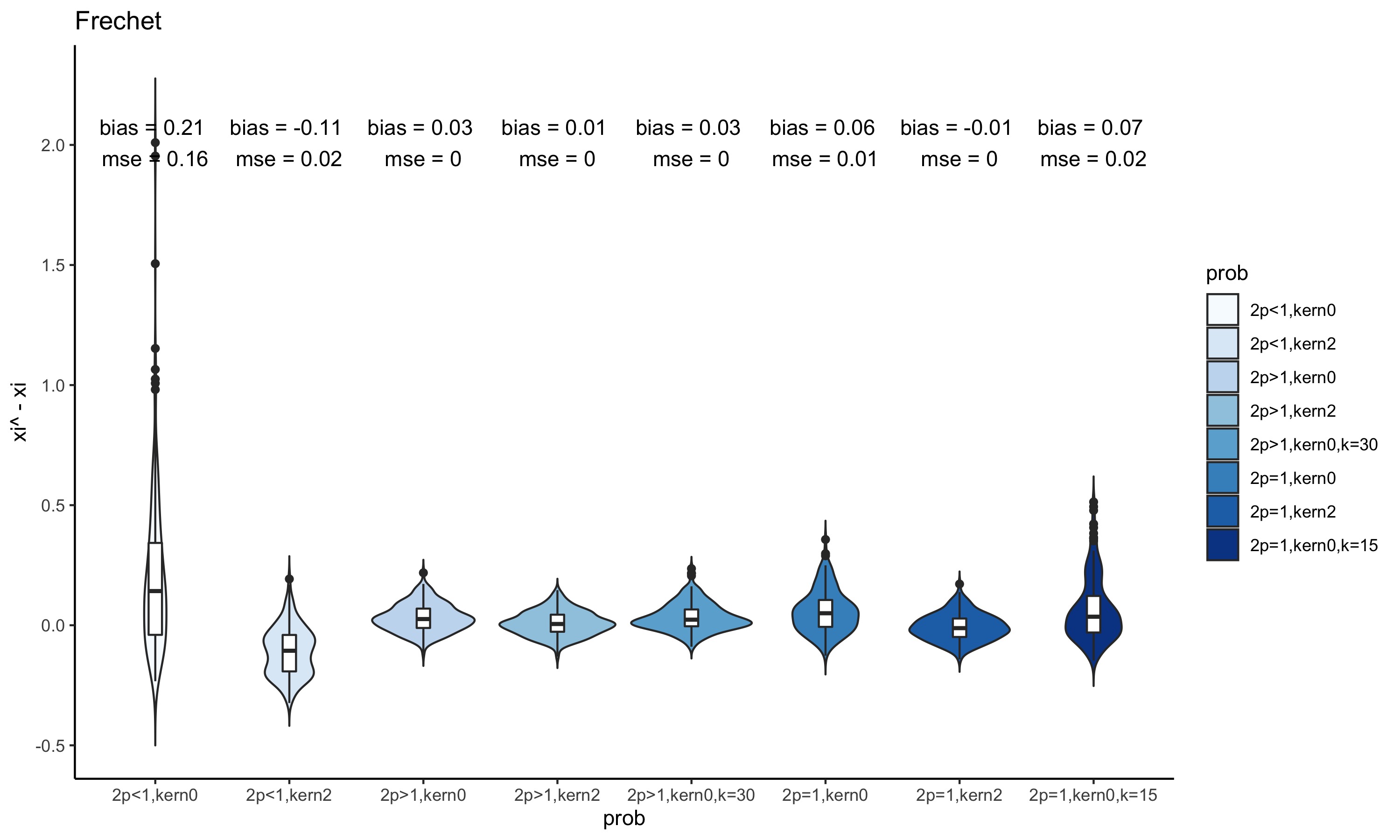}
\caption{\color{black}Violin plots for the simulation results in case of Fr\'echet distributions  under different non-censoring asymptotic probabilities. We take the difference between estimated and theoretical $\xi$. The cases $2p<1$ and $2p\ge 1$ correspond to heavy and light censoring, respectively. The labels kern0 and kern2 correspond to the use of the estimators $H_k^{\mathcal{K}_0}$ and $H_k^{\mathcal{K}_2}$. The specified $k=30,\,15$ are the theoretical optimal sample fractions.}
\label{v2}
\end{figure}

\begin{figure}[hh]
\centering
\includegraphics[width=14cm,trim=0cm 0cm 0cm 0cm,clip]{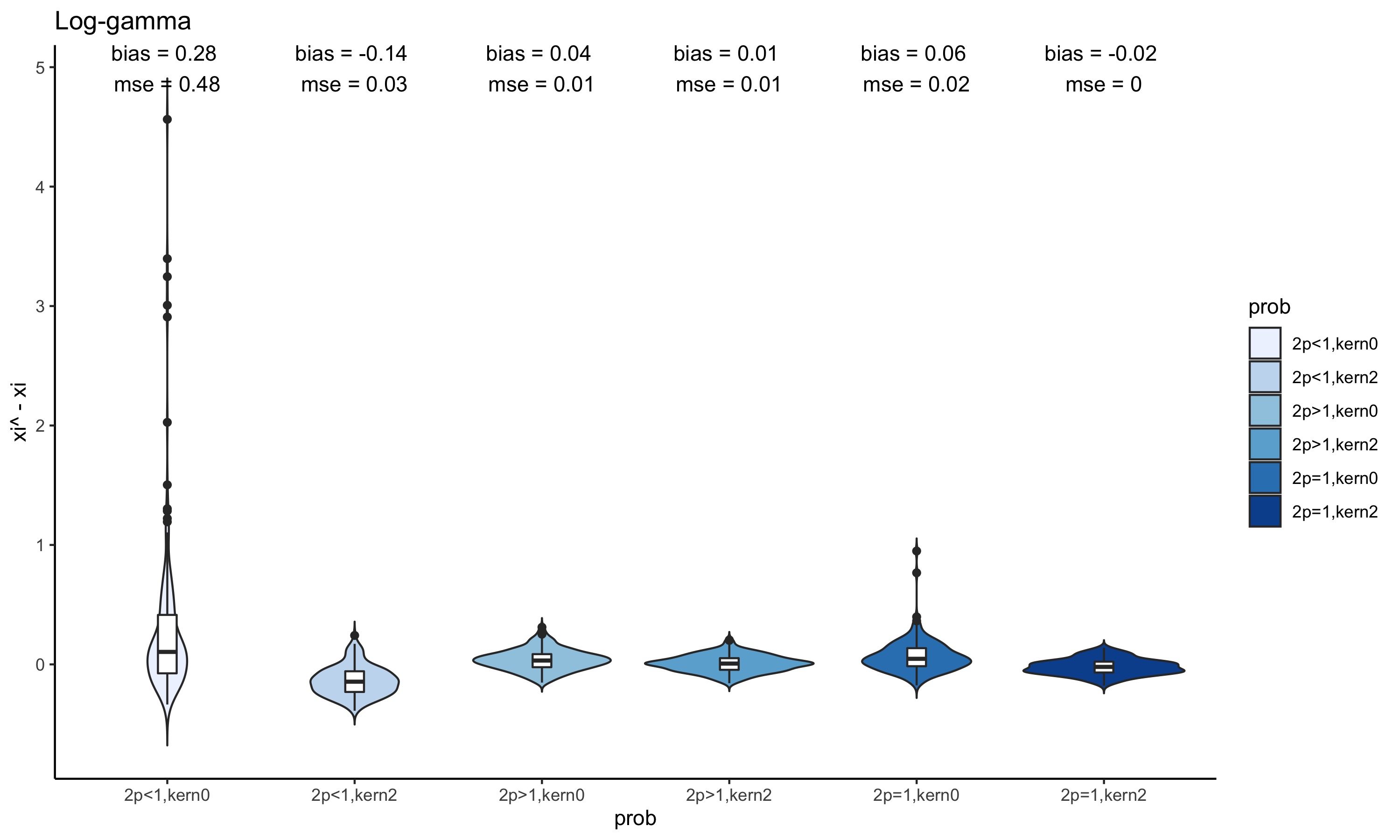}
\caption{\color{black}Violin plots for the simulation results in case of log-gamma distributions  under different non-censoring asymptotic probabilities. We take the difference between estimated and theoretical $\xi$. The cases $2p<1$ and $2p\ge 1$ correspond to heavy and light censoring, respectively. The labels kern0 and kern2 correspond to the use of the estimators $H_k^{\mathcal{K}_0}$ and $H_k^{\mathcal{K}_2}$.}
\label{v3}
\end{figure}

In Figures \ref{v1}, \ref{v2} and \ref{v3} we provide violin plots for the $\mathcal{K}_0$ and $\mathcal{K}_2$-based estimator at the  threshold  selected according to the automatic procedure given in the previous section. We have taken $\rho_z=-1,-3/2,-1/2$ respectively for the three distributions that we consider. These values were permuted (the resulting plots are omitted) and the results were not very sensitive to the choice of $\rho_z$. To avoid degeneracies, a cutoff of $1/5$ of the size of the data set was used for the empirical variance estimates.
 We also add the results of the parameter estimates when taking $k$ fixed at the theoretical optimal value. {\color{black} The latter theoretical optimum is only available in the light-censoring cases and for distributions properly belonging to the Fr\'echet domain of attraction (not the log-gamma case). A general observation is that the violin plots bundle together close to zero, as is desired. Looking closer, we observe that in the regularly varying setup, the adaptive selection of $k$ together with the use of the kernels $\mathcal{K}_0$ and $\mathcal{K}_2$ comes very close to the performance of the Hill estimator evaluated at the theoretical optimum (essentially an oracle estimator, since we input the parameters of the simulated data into this theoretical optimum). As rough guidelines, we observe that the heavy censoring case has significantly worse behaviour than the light censoring case, and that the $\mathcal{K}_2$-based estimator has the best behaviour, agreeing with the conclusions from Figures \ref{bvm_burrs}, \ref{bvm_frechets} and \ref{bvm_LGs}. We believe that this regime-shift between light and heavy censoring is responsible for the increased number of outliers in cases where $p\not > 1/2.$}


\section{Insurance Application: censored claims data vs ultimates}\label{secins}
We now proceed to analyze an insurance dataset consisting of $837$ motor third-party liability (MTPL) insurance claims from $1995$ till $2010$. This data set has been previously described and studied in \cite{abt}, \cite{abbtrim} (without censoring, using ultimate values instead) and \cite{bladt2019combined} (using both censoring and ultimate values). 

\begin{figure}[hh]
\centering
\includegraphics[width=10cm,trim=0cm 0cm 0cm 0cm,clip]{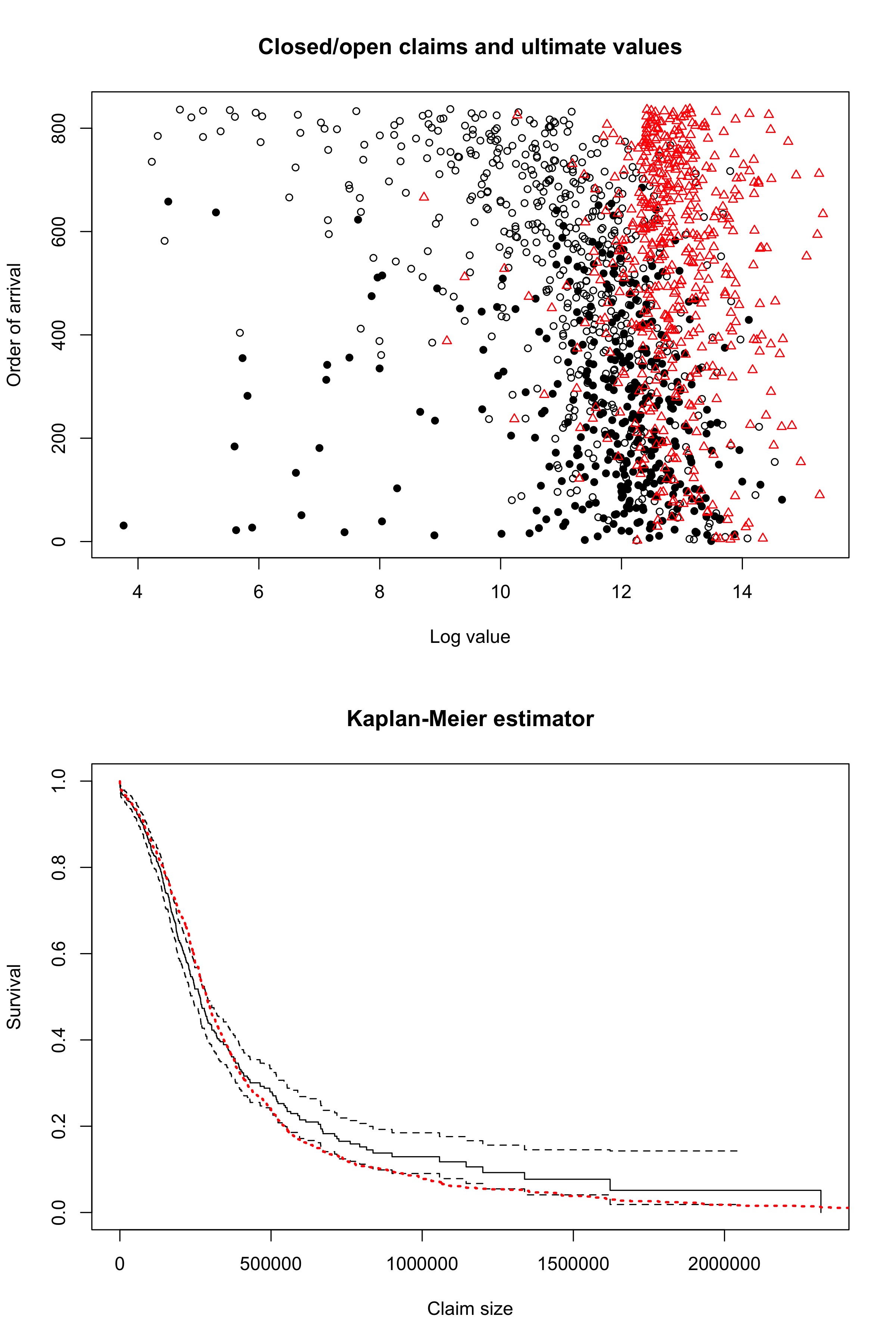}
\caption{MTPL insurance claims data. Top: open (empty circle) and closed (full circle) claims, together with the ultimate values (triangle) for the open claims. Bottom: Kaplan-Meier survival estimate for the open and closed claims, together with the empirical survival function of the ultimates (dotted).}
\label{mtpl_description}
\end{figure}

The data exhibit right-censoring, that is, a claim size is partially observed whenever the development of the claim payment is ongoing and the claim is not yet closed. Closed claim sizes are thus considered as observed data points, and open claims are considered as right-censored observations. In \cite{bladt2019combined} it was argued that the assumption of random censoring and heavy-tailedness is adequate. In the top panel of Figure \ref{mtpl_description} we have the three kinds of data that are available (open claims, closed claims and ultimates), and on the bottom panel the survival function of the open and closed claims using the Kaplan-Meier estimator together with the empirical survival function of the ultimates. We observe that the tail of the latter under-estimates the tail index that is suggested by using survival-analysis techniques.

\begin{figure}[hh]
\centering
\includegraphics[width=11cm,trim=0cm 0cm 0cm 0cm,clip]{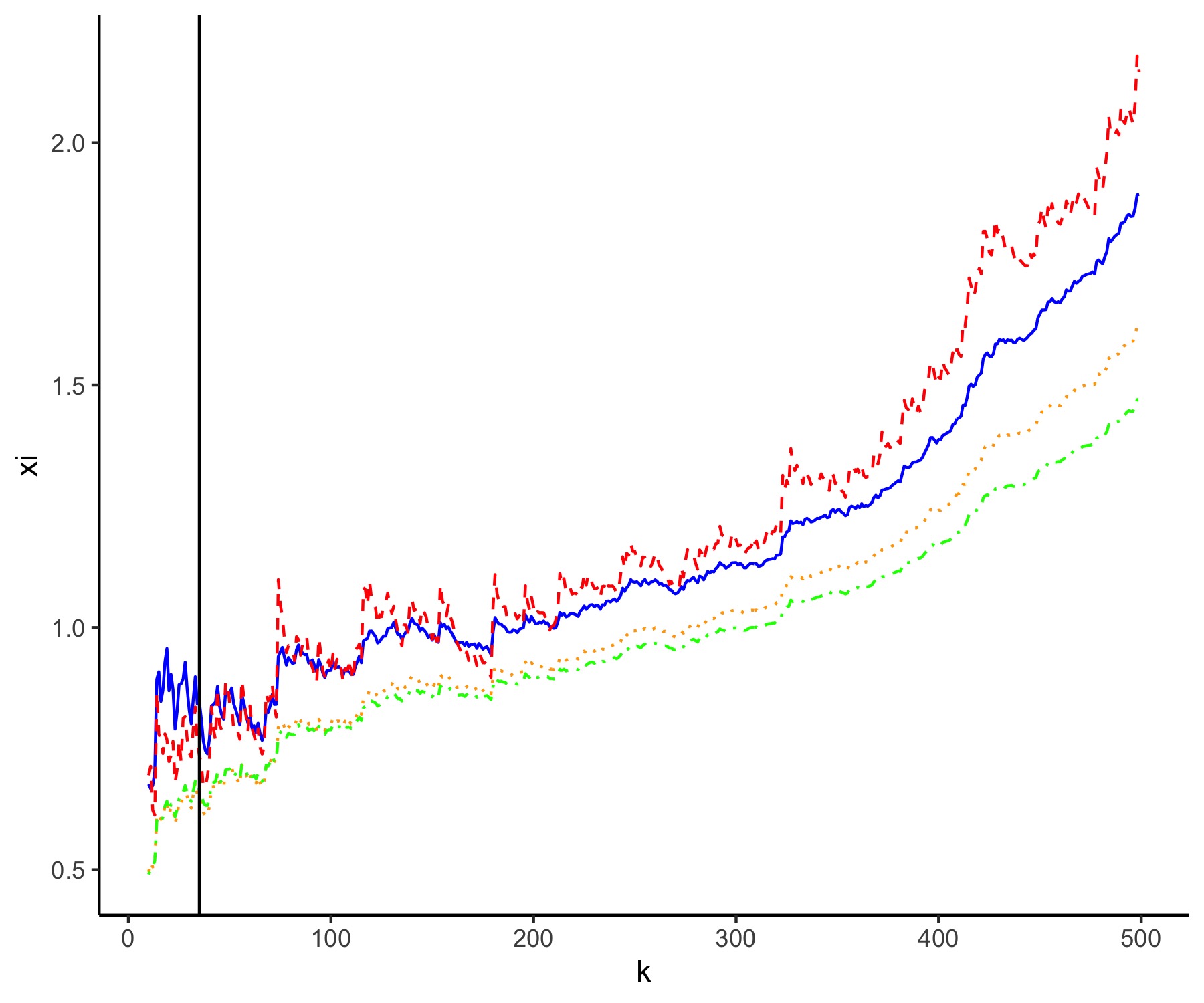}
\caption{Estimates of $\xi$ for the MTPL insurance claim size data: $H_k=H_k^{\mathcal{K}_0}$ in solid blue, $H_k ^A= H_k^{\mathcal{K}_1}$ in dashed red, $H_k^{\mathcal{K}_2}$ in dotted orange and the Worms estimator $H_k^W$  in dashed and dotted green. The vertical line is at the estimated optimal $k$ for $H_k=H_k^{\mathcal{K}_0}$.}
\label{mtpl_xi}
\end{figure}

We now apply the censored tail estimators introduced in this paper to the data. Using the same mechanism as for the simulation study (and $\rho_z=-1$) we find that $k=35$ is the optimal threshold for the estimator using the kernel $\mathcal{K}_0$ (see Figure \ref{mtpl_xi}). As observed in the simulations, it is perfectly reasonable to consider the other estimators ($ H_k^{\mathcal{K}_1}$, $H_k^{\mathcal{K}_2}$ and the Worms estimator $H_k^W$) at this value as well. This yields the estimates
$$\hat\xi=0.841,\:\: 0.761,\:\: 0.654,\:\: 0.664.$$
The latter two values correspond to the kernel $\mathcal{K}_2$ and to the Worms estimator $H_k^W$. Notice that they are quite close, and although the simulations suggest that the two last values are the best performing, the $95\%$ confidence interval for the first of these estimators is given by $[0.415,\, 1.267]$ which amply accommodates all four estimates. Hence, in practice, with only one sample available, it is not possible to make overly conclusive claims regarding the superiority of these point estimates, since they are not statistically distinguishable. The corresponding $99.5\%$ quantile estimators are given by $14957214,\:\:12093195,\:\:9129021,\:\:9355831$, illustrating how small changes in tail estimation can lead to large differences in the quantile scale. Note however that the quantile estimates based on $H_k^W$ and $H_k^{\mathcal{K}_2}$ are quite stable for $k \leq 80$. 

{\color{black}A more refined analysis of $\rho_z$ is known to be unstable, but can be routinely applied (for instance using the $\texttt{mop}$ function from the $\texttt{R}$ package $\texttt{evt0}$). For the Hill estimator of the $Z_i$ variates (ignoring censoring) this gives the estimate $\hat\rho_z=-0.616$, which is relatively close to our choice of $-1$, given the high variability of second-order parameter estimators. Repeating the above analysis with this value has small quantitative influence and no additional qualitative insight, and is thus omitted.}

Previous studies, using the ultimate values (cf. \cite{abbtrim}, with subsequent agreement in \cite{albrecher2019matrix}), that is, internal projected values of the claim sizes at closure, suggested a tail index of about $0.48$. In \cite{bladt2019combined}, combining this expert information with the estimator corresponding to $\mathcal{K}_0$, intermediate values between the purely statistical $0.87$ and the purely expert information $0.48$ were suggested. The present value of $0.66$ is an interesting intermediate value that arises from a purely statistical procedure.

\begin{figure}[hh]
\centering
\includegraphics[width=11cm,trim=0cm 0cm 0cm 0cm,clip]{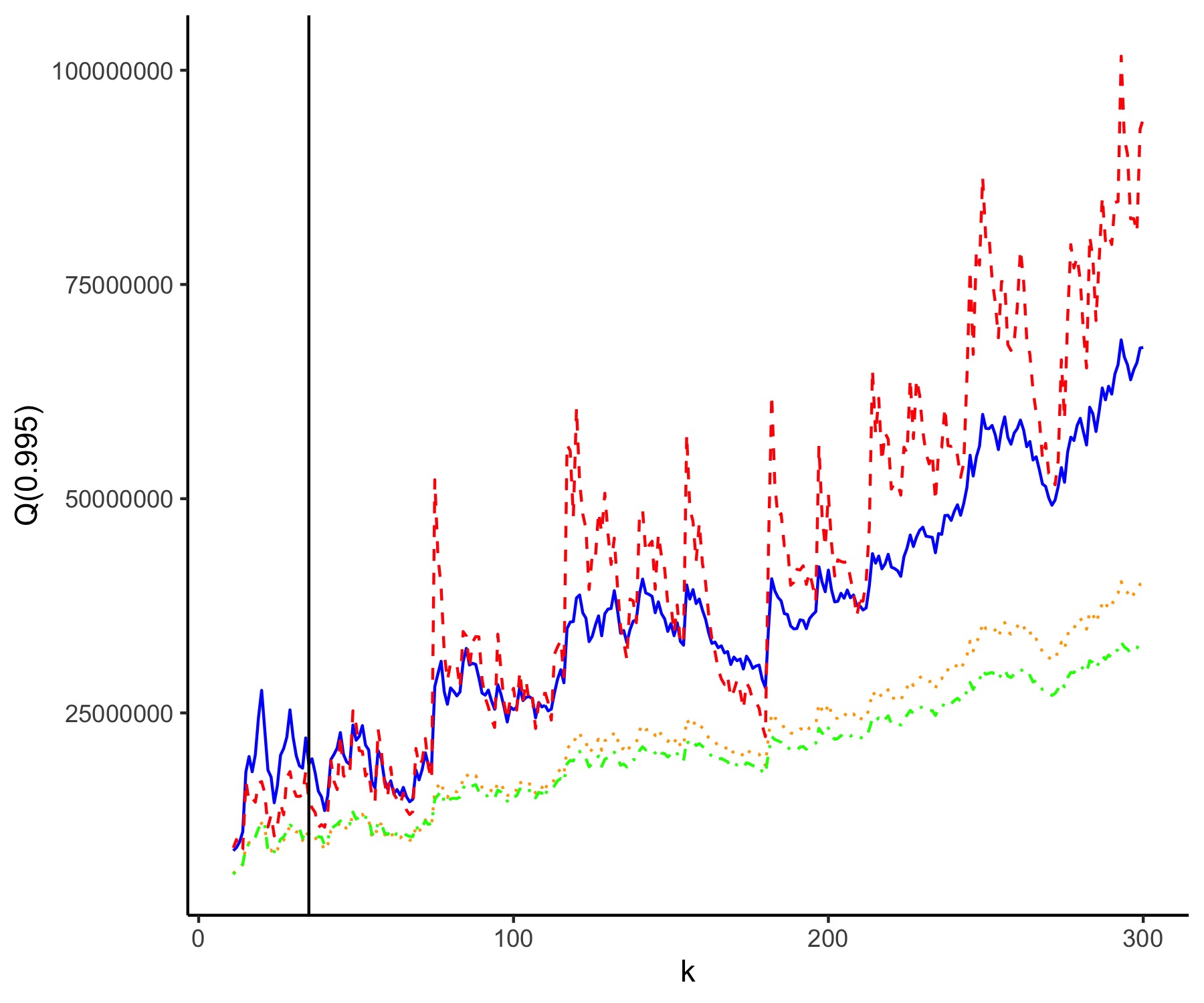}
\caption{Estimates of the $99.5\%$ quantile for the MTPL insurance claim size data based on the tail estimators: $H_k=H_k^{\mathcal{K}_0}$ in solid blue, $H_k ^A= H_k^{\mathcal{K}_1}$ in dashed red, $H_k^{\mathcal{K}_2}$ in dotted orange and the Worms estimator $H_k^W$  in dashed and dotted green. The vertical line is at the estimated optimal $k$ for $H_k=H_k^{\mathcal{K}_0}$.}
\label{mtpl_VaR}
\end{figure}

\section{Conclusion}
In this paper we developed novel extreme value estimators under right-censoring in a kernel framework. The latter class is closed (in the asymptotic sense) under the averaging operation of their trimmed versions, by a simple replacement of kernel. 
The asymptotic behaviour is given for arbitrary kernels, which allows us to compute, for instance, the expression for the MSE as a function of $k$. The choice of the optimal threshold with respect to MSE is explored in connection with the empirical variance of the trimmed trajectories, which leads to an automated way of selecting a threshold. As for the non-censored case, the idea of selecting a threshold by exploiting this link, circumvents the usual estimation difficulties and instabilities which arise in previous approaches in the literature which typically require the estimation of the second-order parameter $D$, and of $\xi$ itself. Despite its simplicity, simulation studies suggest that the method is also efficient. In fact, when compared with the theoretically optimal value, the latter sometimes is too small to be of any practical relevance, and then our adaptive estimator is superior. In the other cases, when the theoretically optimal value is sensible, our estimator also performs well against it. We finally apply the procedure to a well-understood insurance dataset, and the simulation studies suggest that the instances where $\mathcal{K}_0$ has been used in the literature (either alone, or in combination with expert information) to analyze these data could very possibly be improved by considering $\mathcal{K}_2$ instead. Interesting directions for further research include trimming the kernel estimators from above, to remove outliers from data, and to apply combined tail information using censored data and expert information with the new kernels, improving the previous methods. Finally, it will be interesting to consider optimality criteria for the choice of $k$ for any kernel, and to work out criteria for the selection of the optimal kernel from a purely mathematical point of view. \\

\textbf{Acknowledgement.} The authors acknowledge financial support from the Swiss National Science Foundation Project 200021\_191984.

\bibliographystyle{apalike}
\bibliography{Hill_trimming_censored}

\end{document}